# A NOTE ON PARA-QUATERNION MANIFOLDS

Adrian Mihai Ionescu, Gabriel Eduard Vîlcu
*Politehnica University Bucharest, Romania*
*UPG University, Ploieşti, Romania*

Abstract:

We show that a para-quaternion nearly Kahler manifold is necessarily a para-quaternion Kahler manifold

*Introduction*

A para-quaternion Hermitian structure on a pseudo-Riemannian manifold $(M,g)$ is a 3-dimensional subbundle $V$ of $End(TM)$, which locally admits bases of skew-symmetric $I,J,K$ w.r.t. $g$ that further satisfies $I^2 = -J^2 = -K^2 = -Id$, $IJ = -JI = K$, $JK = -KJ = -I$, $KI = -IK = J$.

Let $(M,g,V)$ be an almost para-quaternion-Hermitian manifold (as above), then it is para-quaternion Kahler if $V$ is invariant under covariant differentiation w.r.t. Levi-Civita conexion of $(M,g)$.

Several authors studied para-quaternion geometry and its similitude to the quaternion geometry, see e.g. [1], [2].

It follows in particular that $\dim M \vdots 4$, the metric is neutral, i.e. it has signature $(2m,2m)$, where $\dim M = 4m$ and there is a globally well-defined 4-form $\Omega$ on $M$, given by $\Omega = \Omega_I \wedge \Omega_I - \Omega_J \wedge \Omega_J - \Omega_K \wedge \Omega_K$, where $\Omega_I(X,Y) = g(X,IY)$, $\Omega_J(X,Y) = g(X,JY)$, $\Omega_K(X,Y) = g(X,KY)$

We begin with the following

Lemma 1:

For a (pseudo-)Riemannian manifold $(M,g)$ and a (local) skew-symmetric endomorphism $I$ of $TM$ denote $\Omega_I(X,Y) = g(X,IY)$ and $M_I(X,Y) := (\nabla_X I)Y = \nabla I(Y;X)$. Then:

a) $(\Omega_I \wedge \Omega_I)(X,Y,Z,W) = 2\sum_{Y\ Z\ W} \Omega_I(X,Y)\Omega_I(Z,W)$

$$\nabla_U (\Omega_I \wedge \Omega_I)(X,Y,Z,W) =$$
$$2\sum_{Y \; Z \; W} \{g(X,M_I(U,Y))\Omega_I(Z,W) + \Omega_I(X,Y))g(Z,M_I(U,W))\}, \text{ where the sums}$$

are taken over cyclic permutations of $\{Y,Z,W\}$.

b) $g(M_I(X,Y),Z) + g(M_I(X,Z),Y) = 0$

If $I^2 = \pm Id$, then

c) $M_I(X,IY) = -IM_I(X,Y)$, $g(M_I(X,Y),IY) = g(M_I(X,Y),Y) = 0$
d) $M_I(X,Z) = 0$ for any pair of Hermitian-totally real vectors (i.e. $Z \perp \{X,IX\}$) if and only if $M_I = 0$, i.e. $M$ is Kahler.

If furthermore $I,J,K$ satisfy the para-quaternionic identities, then

e) $M_I(X,JY) = M_K(X,Y) - IM_J(X,Y)$, $M_J(X,IY) = -M_K(X,Y) - JM_I(X,Y)$ (and 4 other similar relations hold)

We can state now the result, see [3]

<u>Theorem 1</u>:

Let $(M,g,V)$ be an almost para-quaternion-Hermitian manifold. Suppose that $M$ is nearly-para-quaternion Kahler, that is $(\nabla_X \Omega)(X,Y,Z,W) = 0$, for any $X,Y,Z,W \in TM$. Then $M$ is para-quaternion Kahler.

Proof:

One must show that for any vectors $X,Y$ it holds that $M_I(X,Y), M_J(X,Y), M_K(X,Y) \in Span\{IY,JY,KY\}$. This will be done by systematic use of formulae in Lemma1.

*Claim1*:

For para-quaternion-Hermitian totally real vectors $X,Y$ (i.e. $X \perp \{Y,IY,JY,KY\}$) one has $M_I(X,Y), M_J(X,Y), M_K(X,Y) \in Q(X) + Q(Y)$.

This is true if $\dim M = 8$; then for any triple $\{X,Y,Z\}$ of pairwise para-quaternion-Hermitian totally real vectors, Lemma1 gives $(\nabla_X \Omega)(X,Y,Z,IX) = -2\|X\|^2 g(Y,M_I(X,Z))$ and similarly for $J,K$. Since $M$ is nearly-para-quaternion Kahler, the Claim follows.

*Claim2*:

For para-quaternion-Hermitian totally real vectors $X,Y$ one has $M_I(X,Y), M_J(X,Y), M_K(X,Y) \perp Q(X)$.

Applying Lemma1 again gives
$$(\nabla_X \Omega)(X,Y,JX,IX) = -2\|X\|^2 g(Y,M_I(X,JX)) - 2\|X\|^2 g(IX,M_J(X,Y))$$
$$+ 2\|X\|^2 g(X,M_K(X,Y)),$$ or after rearranging $(M_K + IM_J - JM_I)(X,Y) \perp X$.
Replacing $Y$ by $KY$ gives again by Lemma1 $(-IM_I + JM_J - 3KM_K)(X,Y) \perp X$.

Similarly one obtains

$$(M_I + KM_J - JM_K)(X,Y) \perp X, (3IM_I + JM_J + KM_K)(X,Y) \perp X,$$

$$(M_J + KM_I - IM_K)(X,Y) \perp X, (-IM_I - 3JM_J + KM_K)(X,Y) \perp X.$$

It follows that $IM_I(X,Y), JM_J(X,Y), KM_K(X,Y) \perp X$, equivalent to the Claim, again by use of Lemma1 and substitutions $Y$ to $IY, JY, KY$.

Now one has $M_I(X,Y), M_J(X,Y), M_K(X,Y) \in Span\{IY, JY, KY\}$, for para-quaternion-Hermitian totally real vectors $X,Y$, since
$$g(M_I(X,Y),Y) = g(M_J(X,Y),Y) = g(M_K(X,Y),Y) = 0.$$

The remaining cases to consider now follows easily since $M_I$-s are tensors while e.g.

$$g(M_I(X,IX), IX) = 0$$

$$g(M_I(X,JX),Y) = -g(M_I(X,Y), JX) = 0$$
$$g(M_I(X,JX), JX) = 0$$
$$g(M_I(X,JX), KX) = g(M_K(X,X) - IM_J(X,X), KX) = g(M_K(X,X), KX) - g(M_J(X,X), JX) = 0$$

This ends the proof of the Theorem.

### References:

[1] E. Garcia-Rio, Y. Mathushita, R.Vazquez-Lorenzo
Paraquaternionic Kahler manifolds, Rocky Jour. Spring 2001

[2] S. Vukmirovic
Para-quaternionic reduction, preprint

[3] A. Swann
Some remarks on quaternion-hermitian manifolds, Archivum Mathematicum (Brno), Tomus 33(1997), 319-351